\documentclass[11pt]{article}
\newcommand{\documentdate}{24 IX 2024}

\usepackage{a4wide,latexsym,graphicx,amsmath,varioref}

\title{Examples of slow convergence for adaptive regularization
  optimization methods are not isolated}

\author{
Philippe L. Toint\thanks{naXys, University of Namur, Namur, Belgium. Email: philippe.toint@unamur.be}
}

\newcommand{\beqn}[1]{\begin{equation}\label{#1}}
\newcommand{\eeqn}{\end{equation}}
\newcommand{\req}[1]{(\ref{#1})}
\newcommand{\ms}{\;\;\;\;}
\newcommand{\tim}[1]{\;\; \mbox{#1} \;\;}

\newcommand{\ii}[1]{\{ 1, \ldots, #1 \}}
\newcommand{\iiz}[1]{\{ 0, \ldots, #1 \}}

\newcommand{\calA}{{\cal A}} 
\newcommand{\calB}{{\cal B}} 
\newcommand{\calO}{{\cal O}} 
\renewcommand{\Re}{\hbox{I\hskip -2pt R}}
\newcommand{\smallRe}{\hbox{\footnotesize I\hskip -2pt R}}
\newcommand{\bigfrac}[2]{\frac{\displaystyle #1}{\displaystyle #2}}
\newcommand{\bigmax}{\displaystyle \max}
\newcommand{\sfrac}[2]{{\scriptstyle \frac{#1}{#2}}}
\newcommand{\half}{\sfrac{1}{2}}
\newcommand{\eqdef}{\stackrel{\rm def}{=}}

\newcounter{algo}[section]
\renewcommand{\thealgo}{\thesection.\arabic{algo}}
\newtheorem{theorem}{Theorem}[section]
\newlength{\thmw}
\newcommand{\algo}[3]{\refstepcounter{algo}
\begin{center}\begin{figure}[htbp]
\framebox[\textwidth]{
\parbox{0.95\textwidth} {\vspace{\topsep}
{\bf Algorithm \thealgo : #2}\label{#1}\\
\vspace*{-\topsep} \mbox{ }\\
{#3} \vspace{\topsep} }}
\end{figure}\end{center}}

\newcommand{\lthm}[2]{\vspace{\baselineskip} 
\noindent\framebox[\textwidth]{\parbox{0.95\textwidth}{
    \begin{theorem} \label{#1} \rm #2 \end{theorem} } } \vspace{\baselineskip} }

\date{\documentdate}

\begin{document}

\maketitle

\begin{abstract}
The adaptive regularization algorithm for unconstrained nonconvex
optimization was shown in \cite{NestPoly06,CartGoulToin11d} to
require, under standard assumptions, at most
$\mathcal{O}(\epsilon^{3/(3-q)})$ evaluations of
the objective function and its derivatives of degrees one and two to
produce an $\epsilon$-approximate critical point of order
$q\in\{1,2\}$. This bound was shown to be sharp in
\cite{CartGoulToin10a,CarmDuchHindSidf20} for $q=1$ and in
\cite{CartGoulToin22} for arbitrary $q \in\{1,2\}$. This note revisits
these results and shows that the example for which slow convergence is
exhibited is not isolated, but that this behaviour occurs for a subset
of univariate functions of nonzero measure.
\end{abstract}

{\small
  \textbf{Keywords:} complexity theory, adaptive regularization methods,
  approximate criticality, nonconvex optimization.
}

\section{Introduction}

Adaptive regularization algorithms for unconstrained nonconvex
minimization have been extensively studied in recent years, focusing
first on the case where a cubically regularized quadratic model is
used (see
\cite{Grie81,NestPoly06,WeisDeufErdm07,GoulRobiThor10,CartGoulToin11d,GoulPorcToin12,
  CartGoulToin12d,CartGoulToin12e,BianLiuzMoriScia15,GrapYuanYuan15a,
  MartRayd17,BergDiouGrat17,BellGuriMori21, CartGoulToin20}, for
instance) and then extending the concept to methods using models of
arbitrary degree \cite{BirgGardMartSantToin17,Duss15}.  We refer the
reader to \cite{CartGoulToin22} for an extensive coverage of this
class of algorithms. A remarkable feature of methods in this class is
their optimal worst-case evaluation complexity: it can indeed be shown
that an $\epsilon$-approximate critical point of order $q$ must be
obtained by an adaptive regularization algorithm using models of
degree $p\geq q$ in at most $\mathcal{O}(\epsilon^{(p+1)/(p-q+1)})$
evaluations of the objective function and its derivatives of degrees 1
to $p$ (see \cite{BirgGardMartSantToin17} for $q=1$ and
\cite{CartGoulToin22} for arbitrary $q\leq p$). This complexity bound
was shown to be sharp in \cite{CartGoulToin10a} for the case where
$p=2$ and $q=1$, and this result was extended in \cite{CartGoulToin22}
to arbitrary $p$ and $q \leq p$ and, independently, in
\cite{CarmDuchHindSidf20} for the special case where $q=1$.  These sharpness
results all hinge on exhibiting of a typically quite contrived
function on which the algorithm under consideration uses exactly as many
evaluations as allowed by the complexity bound to achieve approximate
criticality. Because of this contrived nature, one is naturally led
to the question\footnote{Adressed to the author in several occasions.}
of how ``exceptional'' these examples are. The
purpose of the present note is to show that, albeit not common, such
examples are not isolated, but rather form a set of nonzero measure. To do so,
we focus on the case where $p=2$.

To reach this conclusion, we first need to set the stage in
Section~\ref{defs-s} where we briefly recall the approximate
criticality measures, restate the {\sf AR2} regularization algorithm
and recall the relevant worst-case complexity result. We next revisit,
in Section~\ref{exmpl-s}, the example of slow convergence proposed in
\cite{CartGoulToin10a} and modify it to clarify the amount of
freedom available in its construction. A brief discussion and
conclusion is finally proposed in Section~\ref{concl-s}.

\section{The context}\label{defs-s}

In what follows, we consider the unconstrained minimization problem
\beqn{problem}
\min_{x\in\smallRe^n} f(x)
\eeqn
where $f$ is twice times continously differentiable function from
$\Re^n$ into $\Re$, with globally Lipschitz second derivative, meaning
that there exist a constants $L_j\geq 0$ such that
\[
\|\nabla_x^2 f(x) - \nabla_x^2 f(y) \| \leq L \|x-y\|,
\]
for all $x,y \in \Re^n$. To make the problem well-defined, we also assume that $f$ is
bounded below. We denote by $T_{f,2}(x,s)$ the second
order Taylor's series of such a function computed at $x$ for a step $s$, that is
\[
T_{f,1}(x,s) = f(x)+\nabla_x^1f(x)^Ts
\tim{ and }
T_{f,2}(x,s) = f(x)+\nabla_x^1f(x)^Ts + \half s^T\nabla_x^2f(x)s.
\]
For $j \in \{1,2\}$, we next define
the \textit{$j$-th order criticality measure}
\beqn{phi-def}
\phi_{f,j}(x) = f(x) - \min_{\|d\|\leq 1}T_{f,j}(x,d),
\eeqn
that is the \textit{largest decrease of the $j$-th order Taylor's series of $f$ at $x$
that is achievable in an Euclidean neighbouhood of radius one
centered at $x$}. The measure
$\phi_{f,j}(x)$ is a continuous function of $x$ 
for $j= 1,2$ and we note that it is also independent of $f(x)$. Following
\cite[Section~12.1.5]{CartGoulToin22}, we use this measure to express
the notion of approximate first- and second-order criticality
and say that, given $\epsilon = (\epsilon_1,\ldots,\epsilon_q)$, the point
$x \in \Re^n$ is $\epsilon,$-approximate critical of order $q$ (with
$q \in \{1,2\}$), whenever
\[
\phi_{f,j}(x) \leq \frac{\epsilon_j}{j} \tim{ for } j\in \{1,\ldots,q\}.
\]
Theorem~12.1.3 in \cite{CartGoulToin22}
describes the (natural) way in which $\phi_{f,j}(x)$ relates to
the familiar criticality measures for the values $j=1$ and $q=2$ of
interest here.  In particular,
we have that, for all $x\in\Re^n$,
$\phi_{f,1}(x)= \|\nabla_x^1f(x)\|$ and
  $\phi_{f,2}(x) = | \min[0, \lambda_{\min}[\nabla_x^2f(x)] |$
 whenever $\phi_{f,1}(x)= 0$,
  where $\lambda_{\min}[H]$ denotes the leftmost
  eigenvalue of the symmetric matrix $H$.  
Theorem~12.1.8 in the same reference also gives a lower bound on the values
of (the derivatives of) $f$ in a neighbourhood of $x$.

Having clarified our objective (finding an approximate $q$-order critical
point of $f$), we recall the definition of the {\sf AR$2$} algorithm,
the relevant generalized adaptive regularization minimization method.  As
alluded to in the introduction, this algorithm proceeds, at each iterate $x_k$, by minimizing,
the regularized quadratic model given by
\beqn{model}
m_k(s) = T_{f,2}(x_k,s) + \frac{\sigma_k}{6} |s|^3.
\eeqn
Note that $m_k(s)$ is bounded below, which makes it minimization in
$\Re^n$ well-defined. The minimization of $m_k$ need not be exact
(although we will not use this freedom in what follows) but can be
terminated as soon as the criticality measures $\phi_{m_k,j}$ \textit{associated with
the model} are small enough or the step is large enough (see
\cite{CartGoulToin22} or \cite{GratToin23} for details). 

The {\sf AR2} algorithm (with exact step) is then specified
%\vpageref{AR2-algo}.
as follows.

\algo{AR2-algo}{{\sf AR2} for $\epsilon$-approximate
  $q$-th-order minimizers
}{
\vspace*{-0.3 cm}
\begin{description}
\item[Step 0: Initialization.]
  A criticality order $q\in\{1,2\}$,
  an initial point $x_0$  and an initial regularization parameter $\sigma_0>0$
  are given, as well as accuracy levels  $\epsilon_1\in(0,1]$ and $\epsilon_2\in(0,1]$
  The constants $\theta$, $\eta_1$,
  $\eta_2$, $\gamma_1$, $\gamma_2$, $\gamma_3$ and $\sigma_{\min}$ are also
  given, and satisfy %\vspace*{-2mm}
  \vspace*{-1mm}
  \[
  \begin{array}{c}
    \theta \in (0,1),\;\;
    \sigma_{\min} \in (0, \sigma_0], \;\;
      0 < \eta_1 \leq \eta_2 < 1 %\\
      \tim{and} 0< \gamma_1 < 1 < \gamma_2 < \gamma_3. \vspace*{-1mm}
  \end{array}
  \]
  \vspace*{-1mm}
  Compute $f(x_0)$ and set $k=0$.

\item[Step 1: Test for termination. ]
  For $j = 1,\ldots,q$,
  \begin{enumerate}
  \item Evaluate $\nabla_x^j f(x_k)$ and compute $\phi_{f,j}(x_k)$.
  \item If \vspace*{-2mm}
  \[
  \phi_{f,j}(x_k) \geq  \frac{\epsilon_j}{j},
  \]
  \vspace*{-1mm}
  go to Step~2.
  \end{enumerate}
  Terminate with the approximate solution $x_\epsilon=x_k$.
\item[Step 2: Step calculation. ]
  If not available, evaluate $\nabla^2_x f(x_k)$ and compute a global
  minimizer $s_k$ of the model $m_k(s)$ given by \req{model}.
\item[Step 3: Acceptance of the trial point. ]
  Compute $f(x_k+s_k)$ and define
  \beqn{rhokdef}
  \rho_k = \frac{f(x_k) - f(x_k+ s_k)}{T_{f,2}(x_k,0)- T_{f,2}(x_k,s_k)}.
  \eeqn
%  \vspace*{-1mm}
  If $\rho_k \geq \eta_1$, then define $x_{k+1} = x_k +s_k$.
  Otherwise define $x_{k+1} = x_k$.

\item[Step 4: Regularization parameter update. ]
  Set
  \beqn{sigupdate}
  \sigma_{k+1} \in \left\{ \begin{array}{ll}
  {}[\max(\sigma_{\min}, \gamma_1\sigma_k), \sigma_k ]  & 
    \tim{if} \rho_k \geq \eta_2, \\
  {}[\sigma_k, \gamma_2 \sigma_k ]      &\tim{if} \rho_k \in [\eta_1,\eta_2),\\
  {}[\gamma_2 \sigma_k, \gamma_3 \sigma_k ] & \tim{if} \rho_k < \eta_1.
  \end{array} \right.
  \eeqn
  Increment $k$ by one and go to Step~1 if $\rho_k\geq \eta_1$, or to Step~2 
  otherwise.
\end{description}
}

The evaluation complexity of this algorithm is stated in
\cite[Theorem~12.2.14]{CartGoulToin22}, which we summarize as follows.

\lthm{arqp-complexity}
{
Under the assumptions on $f$ stated at the beginning of this section
and given a criticality order $q\in\{1,2\}$,  the {\sf AR2} algorithm
requires at most  
     \beqn{noncomp-nf-1}
     \calO\left(\bigmax_{j\in\ii{q}}\epsilon_j^{-\frac{3}{3-j}}\right)
     \eeqn
     evaluations of $f$, and its derivatives of orders one and two to produce an
     iterate $x_\epsilon$ such that $\phi_{f,j}(x_\epsilon)\leq \epsilon_j
     /j$ for $j\in\{1,2\}$.
}

\section{The example}\label{exmpl-s}

We now describe a one-dimensional example for which the convergence of
the {\sf AR2} algorithm is ``slow''. Given $q \in \{1,2\}$ and
$\epsilon_1 = \epsilon_2 =\epsilon\in(0,\sfrac{1}{4}]$, we wish to construct a 
function from $\Re$ to $\Re$ such that the {\sf AR2} algorithm, when
applied to minimize this function, will take exactly
\beqn{keps}
k_\epsilon = \left\lceil\epsilon^{-\frac{3}{3-q}}\right\rceil
\eeqn
iterations and evaluations of the objective function and its first and
second derivatives to achieve $\epsilon$-approximate criticality of
order $q$. This function is constructed by first defining sequences of
iterates $\{x_k\}_{k=0}^{k_\epsilon}$, function and derivatives values
$\{\{f_k^{(j)}\}_{j=0}^2\}_{k=0}^{k_\epsilon}$ and associated
``slowly converging'' criticality measures
$\{\{\phi_{f,j}(x_k)\}_{j=1}^2\}_{k=0}^{k_\epsilon}$.  We then
show that these sequences can be viewed as being generated by the {\sf
  AR2} algorithm applied to a piecewise polynomial function from $\Re$ to 
$\Re$ interpolating the values $f_k^{(j)}$.
The general idea is to consider an example in the spirit of
\cite{CartGoulToin10a}, but allowing additional perturbations to the
function's and derivative's values in a way controlled to maintain the
example's mechanism. Note that ensuring termination of the algorithm
at iteration $k_\epsilon$ requires
\beqn{term1-needed}
\phi_{f,1}(x_k) > \epsilon \tim{for} k\in \iiz{k_\epsilon-1}
\tim{and}
\phi_{f,1}(x_{k_\epsilon}) \leq \epsilon.
\eeqn
if $q=1$. If $q=2$, one instead needs that
\beqn{term2-needed-1}
\mbox{either } \phi_{f,1}(x_k) > \epsilon
\tim{ or }
\phi_{f,2}(x_k) \geq \frac{\epsilon}{2} \tim{for} k\in \iiz{k_\epsilon-1}
\eeqn
and 
\beqn{term2-needed-2}
\phi_{f,1}(x_{k_\epsilon}) \leq \epsilon \tim{ and } \phi_{f,2}(x_{k_\epsilon}) \leq \frac{\epsilon}{2}.
\eeqn

We first set $p \eqdef 3-q$ (so that $\{p,q\} = \{1,2\}$)
and define the sequences of derivatives values for $j \in \{1,2\}$
and $k \in \iiz{k_\epsilon}$ by
\beqn{fkj}
f_k^{(j)} =
\left\{ \begin{array}{ll}
  -\beta_{q,k} \alpha_k\epsilon-\alpha_k\epsilon &\tim{for} j = q,\\*[1ex]
  -\beta_{p,k} [\alpha_k\epsilon]^{\frac{q}{p}} &\tim{for} j = p,\\
\end{array}\right..
\eeqn
The sets $\calA  = \{(\alpha_k-1)\}_{k=0}^{k_\epsilon}$ and
$\calB = \{\{\beta_{j,k}\}_{j=1}^2\}_{k=0}^{k_\epsilon}$
define the perturbations of the derivative's values.
That they are acceptable for the example depends crucially on the conditions
\beqn{alpha-conds}
\alpha_k \in [1,2] \tim{for} k\in\iiz{k_\epsilon-1} \tim{and}
  \alpha_{k_\epsilon} = 0
\eeqn
and
\beqn{betaj-cond}
\beta_{q,k} \in[0,\half],\ms
\beta_{p,k} = -\beta_{q,k} \tim{for} k\in\iiz{k_\epsilon-1}
\tim{and}
\beta_{q,k_\epsilon}=\beta_{p,k_\epsilon}=0.
\eeqn
Given \req{fkj}, one verifies that, for sequences $\{x_k\}$ and $\{f_k^{(0)}\}$ yet to be defined,
\beqn{Tfp}
T_{f,2}(x_k,s) %= \sum_{j=0}^2 \frac{f_k^{(j)}}{j!}s^j
= f_k^{(0)} - \alpha_k\epsilon\frac{s^q}{q}
-\beta_{q,k} [\alpha_k\epsilon]\,\frac{s^q}{q}
-\beta_{p,k} [\alpha_k\epsilon]^{\frac{q}{p}}\,\frac{s^p}{p}
\eeqn
so that the model \req{model} is given by
\beqn{mkatsk}
m_k(s) = f_k^{(0)} - \alpha_k\epsilon\frac{s^q}{q}
- \beta_{q,k}[\alpha_k\epsilon]\frac{s^q}{q}
- \beta_{p,k}[\alpha_k\epsilon]^{\frac{q}{p}}\frac{s^p}{p}
+\frac{\sigma_k}{6}|s|^3.
\eeqn
We furthermore choose (as we show below is possible) a constant
regularization parameter
\beqn{sigdef}
\sigma_k = 2 \tim{ for all} k.
\eeqn
Moreover, we consider using the step
\beqn{sk}
s_k =[\alpha_k\epsilon]^{\frac{1}{p}} \tim{for} k\in\iiz{k_\epsilon-1}
\eeqn
and now verify that this step  defines a
global minimizer for the model.
From \req{mkatsk}, \req{sigdef}, \req{sk} and \req{betaj-cond}, we
deduce that
\beqn{zerogradm}
\begin{array}{lcl}
  m_k^{(1)}(s_k)
& = & -\alpha_k\epsilon s_k^{q-1}
  - \beta_{q,k}[\alpha_k\epsilon]s^{q-1}
  - \beta_{p,k}[\alpha_k\epsilon]^{\frac{q}{p}}s^{p-1}
  +\bigfrac{\sigma_k\,s^2}{2}\\
& = & -[\alpha_k\epsilon]^{\sfrac{2}{p}}  
  - \beta_{q,k}[\alpha_k\epsilon]^{\frac{2}{p}}
  - \beta_{p,k}[\alpha_k\epsilon]^{\frac{2}{p}}
  +\bigfrac{\sigma_k\,[\alpha_k\epsilon]^{\sfrac{2}{p}}}{2}\\
& = & -[\alpha_k\epsilon]^{\sfrac{2}{p}}
       \left( 1-\half \sigma_k+\beta_{q,k}+ \beta_{p,k}\right)\\
& = & 0.
\end{array}
\eeqn
and $m_k$ thus admits a global minimizer for $s=s_k$ because the model is the sum of
a quadratic with non-positive slope at $s=0$ and positive cubic
regularization term. We also observe that, because of \req{betaj-cond},
\beqn{relationbeta}
\frac{1}{q}+\frac{\beta_{q,k}}{q} + \frac{\beta_{p,k}}{p}
=\left\{\begin{array}{ll}
  1+\frac{1}{2}\beta_{q,k} \geq 1 & \mbox{if } q=1,\\*[1ex]
  \frac{1}{2} - \frac{1}{2}\beta_{q,k} \geq \frac{1}{4}& \mbox{if } q=2.
\end{array}\right.
\eeqn
As a consequence
\beqn{Tplus}
T_{f,2}(x_k,s_k)
= f_k^{(0)} - [\alpha_k\epsilon]^{\frac{3}{p}}
\left(\frac{1}{q}+ \frac{\beta_{q,k}}{q} + \frac{\beta_{p,k}}{p}\right)
\leq f_k^{(0)} - \frac{s_k^3}{4},
\eeqn
irrespective of the value of $q$.
This gives us the necessary information to define the sequence of
function values by 
\beqn{fk0}
f_0^{(0)} = 3 \times 2^{\frac{3}{p}}
\eeqn
and
\beqn{fplusT}
f_{k+1}^{(0)} = T_{f,2}(x_k,s_k) + \beta_{0,k+1} s_k^3
\eeqn
for some $\beta_{0,k+1}$ (the perturbation on function values) such that
\beqn{beta0-bound}
\beta_{0,0} = 0
\tim{ and }|\beta_{0,k}| \leq \frac{1-\eta_1}{4}
\tim{for} k\in\ii{k_\epsilon}.
\eeqn
As a consequence, we have from \req{rhokdef} and \req{Tplus} that
\[
\left|\rho_k-1\right|
=
\left|\frac{f_{k+1}^{(0)}-T_{f,2}(x_k,s_k)}{f_k^{(0)}-T_{f,2}(x_k,s_k)}\right|
=\leq 4\,|\beta_{0,k+1}| \leq 1-\eta_1,
\]
and thus $\rho_k\geq \eta_1$ for $k\in\iiz{k_\epsilon-1}$, so that every
iteration is successful. Our choice of a constant $\sigma_k$ is
therefore acceptable in view of \req{sigupdate} and $x_{k+1} = x_k+s_k$ for
$k\in\iiz{k_\epsilon-1}$.  We may then set
\beqn{xkdef}
x_0 = 0
\tim{ and }
x_k = \sum_{i=0}^{k-1} s_i
\tim{for}k\in\ii{k_\epsilon}.
\eeqn
A simple calculation shows that
\beqn{phi1}
\phi_{f,1}(x_k) = \left\{\begin{array}{ll}
\alpha_k\epsilon(\beta_{q,k} + 1 ) & \mbox{if } q = 1\\
{}[\alpha_k\epsilon]^2\beta_{p,k}& \mbox{if } q= 2
\end{array}\right.
\eeqn
If $q = 2$, \req{betaj-cond} implies that $\beta_{2,k}+1\geq 1$ and thus
\[
T_{f,2}(x_k,d) = f_k^{(0)}-\beta_{1,k}(\alpha_k\epsilon)^2d-\half \alpha_k\epsilon(\beta_{2,k}+1)d^2
\]
is a concave quadratic in $d$ with nonpositive slope at the origin.
Hence its global minimizer in the interval $[-1,1]$ (see
\req{phi-def}) is achieved for $d=1$, yielding
\beqn{phi2}
\phi_{f,2}(x_k)
= \beta_{1,k}(\alpha_k\epsilon)^2+\half\alpha_k\epsilon(\beta_{2,k}+1)
= (\alpha_k\epsilon)\big(\beta_{p,k}(\alpha_k\epsilon)+\half(\beta_{q,k}+1)\big)
\eeqn
(the value of $\phi_{f,2}(x_k)$ is irrelevant if $q=1$).

\noindent
Let us now consider the values of $\phi_{f,1}(x_k)$ (if $q=1$) or
$\phi_{f,1}(x_k)$ and $\phi_{f,2}(x_k)$ (if $q=2$) as a function of $k$.
Consider the case where $q=1$ first.  Then \req{phi1},
\req{alpha-conds} and \req{betaj-cond} give that
\beqn{term1}
\phi_{f,1}(x_k) > \epsilon \tim{for} k\in \iiz{k_\epsilon-1}
\tim{and}
\phi_{f,1}(x_{k_\epsilon}) = 0
\eeqn
and \req{term1-needed} holds.
For the case $q=2$, \req{betaj-cond}, the
inequality $\epsilon\leq \sfrac{1}{4}$ and \req{alpha-conds} give that
\[
\beta_{p,k}(\alpha_k\epsilon)+\half(\beta_{q,k}+1)
=-\beta_{q,k}(\alpha_k\epsilon)+\half(\beta_{q,k}+1)
\geq \half
\]
Substituting this inequality in \req{phi2} and using \req{alpha-conds}, we obtain that
\beqn{term22}
\phi_{f,2}(x_k) > \half \epsilon \tim{for} k\in \iiz{k_\epsilon-1}
\tim{and}
\phi_{f,2}(x_{k_\epsilon}) = 0,
\eeqn
where we again used \req{alpha-conds} to derive the last equality.
Moreover, from \req{phi1} and \req{alpha-conds},
\beqn{term21}
%\phi_{f,1}(x_k) \leq \epsilon^2 \tim{for} k\in \iiz{k_\epsilon-1}
%\tim{and}
\phi_{f,1}(x_{k_\epsilon}) = 0
\eeqn
(when $q=2$, the value of $\phi_{f,1}(x_k)$ for $k\in
\iiz{k_\epsilon-1}$ is irrelevant since $\phi_{f,2}(x_k) > \half
\epsilon$ for these $k$'s).
Thus termination conditions \req{term2-needed-1} and
\req{term2-needed-2} hold.

\noindent
Turning now to the function values, we observe that \req{relationbeta}
implies that
\[
\frac{1}{q}+ \frac{\beta_{q,k}}{q} + \frac{\beta_{p,k}}{p}
\leq 1 + \left|\frac{1}{q}-\frac{1}{p}\right|\beta_{q,k}
\leq \frac{5}{4},
\]
and therefore, successively using \req{fplusT}, \req{Tplus}, \req{sk},
\req{beta0-bound} and \req{alpha-conds}, that
\[
f_k^{(0)}-f_{k+1}^{(0)}
= [\alpha_k\epsilon]^{\frac{3}{p}}
\left(\frac{1}{q}+ \frac{\beta_{q,k}}{q} + \frac{\beta_{p,k}}{p}\right)
- \beta_{0,k+1}[\alpha_k\epsilon]^{\frac{3}{p}}
\leq [\alpha_k\epsilon]^{\frac{3}{p}}
\left(\frac{5}{4}+\frac{1-\eta_1}{4}\right)
\leq \frac{3}{2}[2\epsilon]^{\frac{3}{p}}.
\]
Thus, using the inequalilty $\epsilon<1$, we obtain that
\beqn{gap-bound}
\begin{array}{lcl}
f_0^{(0)}-f_{k_\epsilon}^{(0)}
& \leq & f_0^{(0)} - \frac{3}{2}\times 2^{\frac{3}{p}}\times k_\epsilon\epsilon^{\frac{3}{p}}\\
& = & f_0^{(0)} - \frac{3}{2}\times 2^{\frac{3}{p}}\times
  \lceil\epsilon^{-\frac{3}{p}}\rceil \epsilon^{\frac{3}{p}}\\
& \leq & f_0^{(0)} - \frac{3}{2}\times 2^{\frac{3}{p}}\times
  (1+\epsilon^{-\frac{3}{p}}) \epsilon^{\frac{3}{p}}\\
& = & f_0^{(0)} - \frac{3}{2}\times 2^{\frac{3}{p}}\times
  (1+\epsilon^{\frac{3}{p}})\\
& < & f_0^{(0)} -3\times 2^{\frac{3}{p}}.
\end{array}
\eeqn
Moreover
\[
f_k^{(0)}-f_{k+1}^{(0)}
= [\alpha_k\epsilon]^{\frac{3}{p}}
\left(\frac{1}{q}+ \frac{\beta_{q,k}}{q} + \frac{\beta_{p,k}}{p}\right)
- \beta_{0,k+1}[\alpha_k\epsilon]^{\frac{3}{p}}
\geq [\alpha_k\epsilon]^{\frac{3}{p}}
\left(\frac{5}{4}-\frac{1-\eta_1}{4}\right)
\geq [\alpha_k\epsilon]^{\frac{3}{p}}
>0,
\]
implying that the sequence $\{f_k^{(0)}\}$ is decreasing.  Thus, in view of
\req{gap-bound} and \req{fk0}, we derive that
\beqn{obj-bound}
f_k^{(0)} \in \left[0,3\times 2^{\frac{3}{p}}\right]
  \tim{for} k\in\iiz{k_\epsilon}.
\eeqn
We now turn to verifying the conditions of
\cite[Theorem~A.9.2]{CartGoulToin22} allowing 
interpolation by a piecewise polynomial with Lipschitz Hessian. We
first verify, using \req{fplusT} and \req{beta0-bound}, that  
\beqn{f0-cond}
| f_{k+1}^{(0)} - T_{f,2}(s_k)| \leq \left(\frac{1-\eta_1}{4} \right)
s_k^3.
\eeqn
We also have from \req{Tfp} that 
\beqn{Tfp1}
\begin{array}{lcl}
T_{f,2}^{(1)}(x_k,s_k) & =
 &-\alpha_k\epsilon s^{q-1}-\beta_{q,k}[\alpha_k\epsilon]s^{q-1} -
\beta_{p,k}[\alpha_k\epsilon]^{\frac{q}{p}}s^{p-1}\\
&=& -[\alpha_k\epsilon]^{\frac{2}{p}}(1+\beta_{q,k}+\beta_{p,k}) = -s_k^2
\end{array}
\eeqn
and
\beqn{Tfp2}
\begin{array}{lcl}
T_{f,2}^{(2)}(x_k,s_k) & =
& -(q-1)\alpha_k\epsilon
s^{q-2}-(q-1)\beta_{q,k}[\alpha_k\epsilon]s^{q-2}
- (p-1)\beta_{p,k}[\alpha_k\epsilon]^{\frac{q}{p}}s^{p-2}\\
& = & -(q-1)[\alpha_k\epsilon]^{\frac{1}{p}}(1+\beta_{q,k})
-(p-1)\beta_{p,k}[\alpha_k\epsilon]^{\frac{1}{p}}\\
&=& s_k[-(q-1)(1+\beta_{q,k})-(p-1)\beta_{p,k}].
\end{array}
\eeqn
Moreover,
\beqn{skratio}
\frac{s_{k+1}^p}{s_k^p}
=\frac{\alpha_{k+1}}{\alpha_k}
\leq 2
\eeqn
because of \req{alpha-conds}.
Using now \req{Tfp1},\req{skratio} and \req{betaj-cond} for $q=1$, we
obtain that, for $k\in \iiz{k_\epsilon-1}$,
\beqn{diff11}
|f_{k+1}^{(1)}-T_{f,2}^{(1)}(x_k,s_k)|
=| -s_{k+1}^2(1+\beta_{q,k}) - s_k^2|
\leq |s_k^2 + 2s_k^2(1+\half)|
= 4\, s_k^2
\eeqn
while using \req{Tfp2},\req{skratio} and \req{betaj-cond} for $q=1$
gives that, for $k\in \iiz{k_\epsilon-1}$,
\beqn{diff21}
|f_{k+1}^{(2)}-T_{f,2}^{(2)}(x_k,s_k)|
=|-\beta_{p,k+1}s_{k+1} + \beta_{p,k}s_k|
\leq \beta_{q,k+1}2^{\sfrac{1}{2}}s_k+\beta_{q,k}s_k
\leq 2\,s_k.
\eeqn
Similarly, using \req{Tfp1},\req{skratio} and \req{betaj-cond} for
$q=2$ yields that, for $k\in \iiz{k_\epsilon-1}$,
\beqn{diff12}
|f_{k+1}^{(1)}-T_{f,2}^{(1)}(x_k,s_k)|
=|-\beta_{p,k+1}s_{k+1}^2-s_k^2|
\leq 4\beta_{q,k+1}s_k^2+s_k^2
%\leq s_k^2(4\beta_{p,k+1}+1)
\leq 3\,s_k^2
\eeqn
and using \req{Tfp2},\req{skratio} and \req{betaj-cond} for
$q=2$ implies that, for $k\in \iiz{k_\epsilon-1}$,
\beqn{diff22}
\begin{array}{lcl}
|f_{k+1}^{(2)}-T_{f,2}^{(2)}(x_k,s_k)|
& =& |-s_{k+1}(1+\beta_{q,k+1})+s_k(1+\beta_{q,k})|\\
& \leq & s_k[ 2(1+\beta_{q,k+1}) + 1+\beta_{q,k}]\\
& < & \sfrac{9}{2}\,s_k
\end{array}
\eeqn
Combining \req{diff11}--\req{diff22}, we see that, in all cases and
for all $k\in \iiz{k_\epsilon-1}$, 
\beqn{diffall}
|f_{k+1}^{(1)}-T_{f,2}^{(1)}(x_k,s_k)| \leq \sfrac{9}{2}\,s_k^2
\tim{and}
|f_{k+1}^{(2)}-T_{f,2}^{(2)}(x_k,s_k)| \leq \sfrac{9}{2}\,s_k.
\eeqn
In addition, \req{fkj}, \req{alpha-conds}, \req{betaj-cond}, the fact that
$\epsilon\leq \half$ and \req{obj-bound} yield that, for all $j\in \{0,1,2\}$ and
all $k\in\iiz{k_\epsilon}$,
\beqn{fkj-bound}
|f_k^{(j)}| \leq \max\left[\,\frac{5}{2}, \,3 \times 2^{\frac{3}{p}}\,\right]
\tim{ and }
|s_k| \leq 1,
\eeqn
where we used \req{sk} and \req{alpha-conds} to derive the second inequality.

\noindent
Combining \req{f0-cond}, \req{diffall} and \req{fkj-bound} , we may
now apply \cite[Theorem~A.9.2]{CartGoulToin22} for Hermite
interpolation with 
\[
\kappa_f = \max\left[ \,\frac{9}{2}, \,\frac{5}{2}, \,3 \times 2^{\frac{3}{p}}, \,1 \,\right]
= 9\times2^{\frac{3}{p}-1}
\approx\left\{\begin{array}{ll}
8.485 & \mbox{if } q =1\\
24 & \mbox{if } q =2
\end{array}\right.
\]
and deduce the existence of a twice times continuously differentiable
piecewise polynomial function $f_{\calA,\calB}$ from $\Re$
to $\Re$ with bounded continuous derivatives of degrees zero to two and
Lipschitz continuous Hessian, which interpolates the data given by
$\{\{f_k^{(j)}\}_{j=0}^2\}_{k=0}^{k_\epsilon}$ at the iterates
$\{x_k\}_{k=0}^{k_\epsilon}$. Moreover, its Hessian's Lipschitz
constant $L$ only depends on $\kappa_f$. The same theorem, together
with \req{obj-bound}, also ensures that  $|f_{\calA,\calB}^{(j)}(x)|$
is uniformly bounded for all $j\in \{0,1,2\}$ and all $x\in
\Re$. Finally, \req{sk}, \req{alpha-conds}, the  inequality $\epsilon
\leq \frac{1}{4}$ and \req{xkdef} imply that
\beqn{xkinter}
x_k \in [0, \half k_\epsilon] \tim{ for all } k \in \iiz{k_\epsilon},
\eeqn
irrespective of the choice of $\calA$ and $\calB$.
This conclude the construction of our example.

\section{Discussion}\label{concl-s}

We note that the example of \cite[Section~12.2.2.4]{CartGoulToin22}
corresponds to selecting
\[
\beta_{k,j} = 0
\tim{ and  }
\alpha_k = 1+\frac{k_\epsilon-k}{k_\epsilon}
\]
for all $k$ and $j$ in the above development\footnote{And also corrects a minor error in equation
(12.2.78) of \cite{CartGoulToin22}.}. 

However, the very fact that we may choose
\beqn{ranges}
\beta_{0,k+1}\in \left[-\frac{1-\eta_1}{4},\frac{1-\eta_1}{4}\right],
\ms
\beta_{p,k}=-\beta_{q,k} \in [0, \half]
\tim{ and }
\alpha_k\in[1,2]
\eeqn
for each $k$ and each $j$ (see \req{alpha-conds} \req{betaj-cond} and
\req{beta0-bound}) tells us, in conjunction with \req{fkj} and
\req{fplusT}, that the interpolation data at $x_k$ may be chosen, for
$k\in\iiz{k_\epsilon-1}$, arbitrarily in an interval of radius
$2(1-\eta_1)s_k^{p+1}/4$ for the objective function see
\req{f0-cond}), and in intervals of radius at least $\epsilon$ for the
derivative of degree $q$ and $\half\epsilon^{q/p}$ for the derivative
of degree $p$ (see \req{fkj}--\req{betaj-cond}). Moreover, the proof
of \cite[Theorem~A.9.2]{CartGoulToin22} reveals that, in each interval
$[x_k,x_{k+1}]$, $f_{\calA,\calB}$ is a linear combination of the
monomials $1$, $s$, $s^2$, $s^3$, $s^4$ and $s^5$ whose coefficients
depend continuously and bijectively on the values of $f_k^{(j)}$ and
$f_{k+1}^{(j)}$ for $j\in\{0,1,2\}$.
Now observe that, in our example,
\[
\left( \begin{array}{c}
  f_k^{(0)}\\f_k^{(q)}\\f_k^{(p)}
\end{array}\right)
= \left(\begin{array}{c}
  T_{2,f}(x_{k-1},s_{k-1})+\beta_{0,k}[\alpha_{k-1}\epsilon]^{\sfrac{3}{p}}\\
  -[\alpha_k\epsilon](1+\beta_{q,k})\\
  {}[\alpha_k\epsilon]^{\sfrac{q}{p}}\beta_{q,k}
\end{array}\right)
\eqdef \Theta_k(\beta_{0,k}, \alpha_k, \beta_{q,k}).
\]
Since the determinant of the Jacobian of $\Theta_k$ is
given by
\[
\left|\left(\begin{array}{ccc}
  [\alpha_{k-1}\epsilon]^{\sfrac{3}{p}} & 0 & 0 \\
  0 & -\epsilon(1+\beta_{q,k}) &  -\alpha_k\epsilon\\
  0 & \bigfrac{q}{p}\beta_{q,k}\epsilon[\alpha_k\epsilon]^{\sfrac{q}{p}-1} & [\alpha_k\epsilon]^{\sfrac{q}{p}}\\
\end{array}\right)\right|
= - [\alpha_{k-1}\epsilon]^{\sfrac{3}{p}}\epsilon^{\sfrac{3}{p}}\alpha_k^{\frac{q}{p}}\left[1+\left(1-\frac{q}{p}\right)\beta_{q,k}\right]
\]
and is nonzero for all values of $(\beta_{0,k}, \alpha_k,
\beta_{q,k})$ in the ranges \req{ranges}, $J_k$ is nonsingular in these ranges and we deduce that
$\Theta_k$ is continuous and bijective for all $k$.
As a consequence, different values of $\beta_{0,k}$, $\alpha_k$ and
$\beta_{q,k}$ for a at least one $k$ result in piecewise polynomials $f_{\calA,\calB}$ which
(continuoulsy) differ in the neighbourhood (in $x$) of at least one
iterate, the radius of this neighbourhood depending on the
parameter-independent Lipschitz constant $L$. We may therefore
conclude that

\vspace*{2mm}
\noindent
\centerline{\fbox{\parbox{0.98\linewidth}{
     \emph{the set of functions on which the {\sf AR2} algorithm can
     take $\calO(\epsilon^{-3/p})$ evaluations is of nonzero measure
     (in the standard topology for continuous functions from $\Re$ to $\Re$)}
}}}

\vspace*{2mm}
\noindent
although the measure of this set may shrink when $\epsilon$ tends to zero.

\vspace*{2mm}
\begin{figure}[htb] % produced by examples.m
\centerline{\includegraphics[height=6cm,width=12cm]{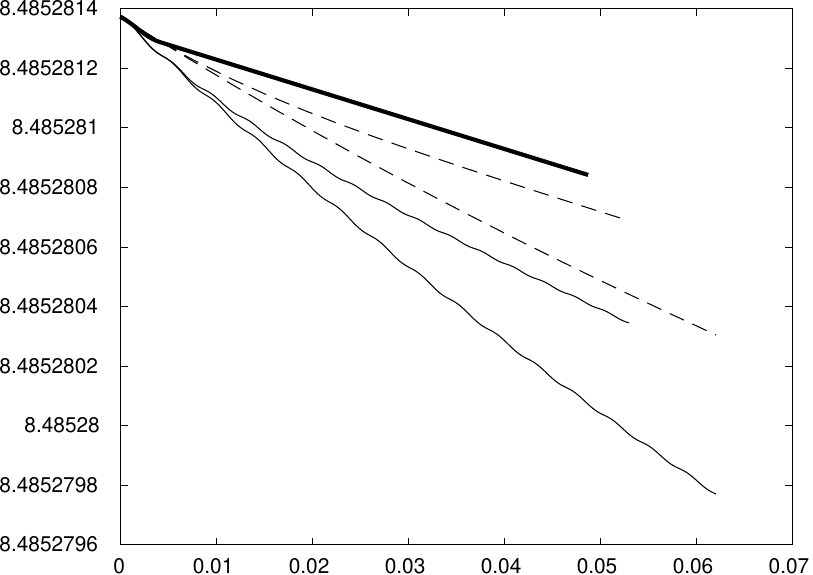}}
\caption{\label{figure:slowex}A few members of the set of functions causing slow
  convergence of the adaptive regularization algorithm ($\epsilon =
  10^{-5}$, $q=1$, showing the first 15 iterations)}
\end{figure}

Figure~\ref{figure:slowex} illustrates the diversity of possible
examples by showing the graphs of a few $f_{\calA,\calB}$
for $q=1$ corresponding to different choices of $\calA$ and
$\calB$ with $\beta_{0,k}=0$ for all $k$. The unpertubed example (i.e.\ with $\alpha_k=1$ and
$\beta_{j,k}=0$ for all $k$ is the top fatter curve while the bottom
dashed one corresponds to setting $\alpha_k=2$ and $\beta_{q,k}=0$ for
all $k$. The thin continuous curves show what happens if $\beta_{q,k}=
\half$ is used for all $k$ with $\alpha_k=1$ or $\alpha_k=2$. Curves have different maximum
abscissa for the same number of iterations because the steplength
$s_k$ given by \req{sk} varies with $\alpha_k$.

The reader may also wonder if our example contradicts the result of
\cite{GratSimToin24}, which states that the {\sf AR2} algorithm
terminates in fact in $o(\epsilon^{3/p})$ rather that in
$O(\epsilon^{3/p})$ iterations and evaluations.  Fortunately, this is
not the case, because the contexts in which these results are derived
differ. More specifically, the example presented above constructs a
different function for every choice of $\epsilon$ (and possible
perturbations).  It states that, if $\epsilon$ is fixed, then there
exists a set of sufficiently smooth functions depending on $\epsilon$
causing slow convergence of the {\sf AR2} algorithm. In contrast, the
result of \cite{GratSimToin24} states that, for any sufficiently
smooth (fixed) function, the number of iterations (and evaluations)
required to achieve $q$-th order $\epsilon$-approximate criticality
\emph{for this function} increases slower than $\epsilon^{-3/p}$ when
$\epsilon$ tends to zero.

\section*{\footnotesize Acknowledgement}

{\footnotesize
Thanks to Geovani Grapiglia for persuading the author to write this
note during an interesting discussion at the Brazopt XIV conference in
Rio de Janeiro, and to the organizers of this meeting for their
support. Thanks also to Serge Gratton for his encouragements, and to
Benedetta Morini, Stefania Bellavia and Margherita Porcelli for
helping to check the details.

%\bibliography{/home/philippe/bibs/refs}

\begin{thebibliography}{10}

\bibitem{BellGuriMori21}
S.~Bellavia, G.~Gurioli, and B.~Morini.
\newblock Adaptive cubic regularization methods with dynamic inexact {H}essian
  information and applications to finite-sum minimization.
\newblock {\em IMA Journal of Numerical Analysis}, 41(1):764--799, 2021.

\bibitem{BergDiouGrat17}
E.~Bergou, Y.~Diouane, and S.~Gratton.
\newblock On the use of the energy norm in trust-region and adaptive cubic
  regularization subproblems.
\newblock {\em Computational Optimization and Applications}, 68:533--554, 2017.

\bibitem{BianLiuzMoriScia15}
T.~Bianconcini, G.~Liuzzi, B.~Morini, and M.~Sciandrone.
\newblock On the use of iterative methods in cubic regularization for
  unconstrained optimization.
\newblock {\em Computational Optimization and Applications}, 60(1):35--57,
  2015.

\bibitem{BirgGardMartSantToin17}
E.~G. Birgin, J.~L. Gardenghi, J.~M. Mart\'{i}nez, S.~A. Santos, and Ph.~L.
  Toint.
\newblock Worst-case evaluation complexity for unconstrained nonlinear
  optimization using high-order regularized models.
\newblock {\em Mathematical Programming, Series~A}, 163(1):359--368, 2017.

\bibitem{CarmDuchHindSidf20}
Y.~Carmon, J.~C. Duchi, O.~Hinder, and A.~Sidford.
\newblock Lower bounds for finding stationary points {I}.
\newblock {\em Mathematical Programming, Series~A}, 184:71--120, 2020.

\bibitem{CartGoulToin10a}
C.~Cartis, N.~I.~M. Gould, and Ph.~L. Toint.
\newblock On the complexity of steepest descent, {N}ewton's and regularized
  {N}ewton's methods for nonconvex unconstrained optimization.
\newblock {\em SIAM Journal on Optimization}, 20(6):2833--2852, 2010.

\bibitem{CartGoulToin11d}
C.~Cartis, N.~I.~M. Gould, and Ph.~L. Toint.
\newblock Adaptive cubic overestimation methods for unconstrained optimization.
  {P}art {II}: worst-case function-evaluation complexity.
\newblock {\em Mathematical Programming, Series~A}, 130(2):295--319, 2011.

\bibitem{CartGoulToin12d}
C.~Cartis, N.~I.~M. Gould, and Ph.~L. Toint.
\newblock Complexity bounds for second-order optimality in unconstrained
  optimization.
\newblock {\em Journal of Complexity}, 28:93--108, 2012.

\bibitem{CartGoulToin12e}
C.~Cartis, N.~I.~M. Gould, and Ph.~L. Toint.
\newblock On the evaluation complexity of cubic regularization methods for
  potentially rank-deficient nonlinear least-squares problems and its relevance
  to constrained nonlinear optimization.
\newblock {\em SIAM Journal on Optimization}, 23(3):1553--1574, 2013.

\bibitem{CartGoulToin20}
C.~Cartis, N.~I.~M. Gould, and Ph.~L. Toint.
\newblock A concise second-order evaluation complexity for unconstrained
  nonlinear optimization using high-order regularized models.
\newblock {\em Optimization Methods and Software}, 35(2):243--256, 2020.

\bibitem{CartGoulToin22}
C.~Cartis, N.~I.~M. Gould, and Ph.~L. Toint.
\newblock {\em Evaluation complexity of algorithms for nonconvex optimization}.
\newblock Number~30 in MOS-SIAM Series on Optimization. SIAM, Philadelphia,
  USA, June 2022.

\bibitem{Duss15}
J.~P. Dussault.
\newblock Simple unified convergence proofs for the trust-region and a new
  {ARC} variant.
\newblock Technical report, University of Sherbrooke, Sherbrooke, Canada, 2015.

\bibitem{GoulPorcToin12}
N.~I.~M. Gould, M.~Porcelli, and Ph.~L. Toint.
\newblock Updating the regularization parameter in the adaptive cubic
  regularization algorithm.
\newblock {\em Computational Optimization and Applications}, 53(1):1--22, 2012.

\bibitem{GoulRobiThor10}
N.~I.~M. Gould, D.~P. Robinson, and H.~S. Thorne.
\newblock On solving trust-region and other regularised subproblems in
  optimization.
\newblock {\em Mathematical Programming, Series~C}, 2(1):21--57, 2010.

\bibitem{GrapYuanYuan15a}
G.~N. Grapiglia, J.~Yuan, and Y.~Yuan.
\newblock On the convergence and worst-case complexity of trust-region and
  regularization methods for unconstrained optimization.
\newblock {\em Mathematical Programming, Series~A}, 152:491--520, 2015.

\bibitem{GratSimToin24}
S.~Gratton, C.-K. Sim, and Ph.~L. Toint.
\newblock Refining asymptotic complexity bounds for nonconvex optimization
  methods, including why steepest descent is $o(\epsilon^{-2})$ rather than
  $\mathcal{O}(\epsilon^{-2})$.
\newblock arXiv:2408.09124, 2024.

\bibitem{GratToin23}
S.~Gratton and Ph.~L. Toint.
\newblock Adaptive regularization minimization algorithms with non-smooth
  norms.
\newblock {\em IMA Journal of Numerical Analysis}, 43(2):920--949, 2023.

\bibitem{Grie81}
A.~Griewank.
\newblock The modification of {N}ewton's method for unconstrained optimization
  by bounding cubic terms.
\newblock Technical Report NA/12, Department of Applied Mathematics and
  Theoretical Physics, University of Cambridge, Cambridge, United Kingdom,
  1981.

\bibitem{MartRayd17}
J.~M. Mart\'{\i}nez and M.~Raydan.
\newblock Cubic-regularization counterpart of a variable-norm trust-region
  method for unconstrained minimization.
\newblock {\em Journal of Global Optimization}, 68:367–385, 2017.

\bibitem{NestPoly06}
{Yu}. Nesterov and B.~T. Polyak.
\newblock Cubic regularization of {N}ewton method and its global performance.
\newblock {\em Mathematical Programming, Series~A}, 108(1):177--205, 2006.

\bibitem{WeisDeufErdm07}
M.~Weiser, P.~Deuflhard, and B.~Erdmann.
\newblock Affine conjugate adaptive {N}ewton methods for nonlinear
  elastomechanics.
\newblock {\em Optimization Methods and Software}, 22(3):413--431, 2007.

\end{thebibliography}
%\bibliographystyle{plain}

}
\end{document}